\documentclass[reqno, 12pt]{amsart}

\usepackage{amsfonts}
\usepackage{amssymb,amsmath}
\usepackage{amsthm}
\usepackage{upref}
\usepackage{enumerate}
\usepackage{bbm}

\makeatletter
\def\LaTeX{\leavevmode L\raise.42ex
    \hbox{\kern-.3em\size{\sf@size}{0pt}\selectfont A}\kern-.15em\TeX}
 \makeatother

\numberwithin{equation}{section}

\newtheorem{lemma}{Lemma}[section]
\newtheorem{theorem}[lemma]{Theorem} 
\newtheorem{corollary}[lemma]{Corollary}
\newtheorem{proposition}[lemma]{Proposition}

\theoremstyle{definition}

\newtheorem{definition}[lemma]{Definition}
\newtheorem{example}[lemma]{Example}

\newtheorem{remark}[lemma]{Remark}

\newcommand{\1}{\mathbbm{1}}

  \newcommand{\e}{\eqref}

\newcommand{\q}{\quad}

\newcommand{\wt}{\widetilde}
\newcommand{\wh}{\widehat}
\newcommand{\la}{\langle}
\newcommand{\ra}{\rangle}
\newcommand{\z}{\zeta}
\newcommand{\ov}{\overline}
 \renewcommand{\d}{\delta}
 
 \newcommand{\cl}{\operatorname{clos}}

 \newcommand{\Ker}{\operatorname{Ker}}

\renewcommand\Im{\operatorname{Im}}
\renewcommand\Re{\operatorname{Re}}

\newenvironment{pf}{\begin{proof}}{\end{proof}}

\def\qqq{\mathrel{\subset\mkern-15mu\lower.38ex\hbox{${\scriptscriptstyle\rightarrow}$}}}

\let\cal\mathcal

\let\Bbb\mathbb

\begin{document}

\title {On semibounded Wiener-Hopf operators}

\author{ D. R. Yafaev}
\address{ IRMAR, Universit\'{e} de Rennes I\\ Campus de
  Beaulieu, 35042 Rennes Cedex, FRANCE}
\email{yafaev@univ-rennes1.fr}

\begin{abstract} 
 We show that a semibounded Wiener-Hopf  quadratic form   is closable in the space $L^2({\Bbb R}_{+})$ if and only if
 its integral  kernel is the Fourier transform of an absolutely continuous measure.
 This allows us to define semibounded Wiener-Hopf  operators and their symbols  under minimal assumptions on their integral  kernels. 
 Our proof relies on   a continuous analogue of the Riesz Brothers theorem   obtained in the paper.
     \end{abstract}

\maketitle


\section{Introduction. Main results }  

{\bf 1.1.}
Wiener-Hopf operators   $W$ can formally be defined in the space $L^2({\Bbb R}_{+})$ of functions  $f(x)$
 by the formula
  \begin{equation}
(W f) (x)=\int_{ {\Bbb R}_{+}}   w(x-y) f(y) dy.
 \label{eq:HD}\end{equation}
 Thus the integral   kernel $w$ of a  Wiener-Hopf   operator depends  on the difference of the variables $x$ and $y$ only. So it is natural to expect that properties of  Wiener-Hopf  operators are close to those of convolution operators acting in the space $L^2({\Bbb R})$.
 
To be precise, we define the operator $W$ via its quadratic form
   \begin{equation}
w[f,f] =\int_{{\Bbb R}_{+}}  \int_{{\Bbb R}_{+}}     w(x-y) f(y) \ov{f (x)}dx dy .
 \label{eq:QFq}\end{equation}
  With respect to   $w$, we a priori only assume that it is a distribution in the class $  C_{0}^\infty ({\Bbb R})'$ dual to $  C_{0}^\infty ({\Bbb R})$. Then the quadratic form is correctly defined for all  $f\in C_{0}^\infty ({\Bbb R}_{+})$.
  This is discussed in more detail in Subsection~2.2.


 To a large extent, the theory of  Wiener-Hopf  operators is parallel to the theory of  Toeplitz operators $T$ acting in the space $\ell^2({\Bbb Z}_{+})$ of sequences by the formula
   \begin{equation}
(T g)_{n}=\sum_{m= 0}^\infty  t_{n-m} g_{m}.
 \label{eq:HDT}\end{equation}
 Roughly speaking, $W$ (also sometimes called Toeplitz operators) and $T$ are continuous and discrete versions of the same object (see Subsection~6.1). To avoid confusion, we use the terms ``Wiener-Hopf$\,$" and  ``Toeplitz" for the operators defined by formulas \e{eq:HD} and  \e{eq:HDT}, respectively. 
 
 However optimal results on 
 Wiener-Hopf  operators are not direct consequences of the corresponding results for Toeplitz operators and, in some sense, they are more general.  One of the differences is that  Wiener-Hopf  operators  require a consistent work with distributions.   As an example, let us state a necessary and sufficient condition for a Wiener-Hopf operator  $W$ to be bounded. 
Note  that  {\it the Fourier transform is always understood in terms of the Schwartz space} ${\cal S}'  $ dual to the space
${\cal S} ={\cal S} ({\Bbb R})$ of   rapidly decaying $C^{\infty}$ functions.
 
  \begin{theorem}\label{B-H} 
  Let the form $w[f,f] $ be defined by the relation \e{eq:QFq} where the distribution $w\in C_{0}^\infty({\Bbb R})'$. Then the estimate
   \begin{equation}
|w[f,f] | \leq C \| f\|^2_{ L^2 ({\Bbb R}_{+} )} , \q \forall f \in C_{0}^\infty({\Bbb R}_{+}),
 \label{eq:QFq1}\end{equation}
 with some  constant $C>0$
 is equivalent to the representation
    \begin{equation}
 w(x) = \frac{1}{2\pi}\int_{ \Bbb R} e^{-i\lambda x} \varphi (\lambda) d\lambda  \q {\rm where}  \q \varphi\in L^\infty ({\Bbb R} ).
 \label{eq:BH}\end{equation}
 Moreover,
  \begin{equation}
\sup_{\| f\|_{ L^2 ({\Bbb R}_{+} )}=1 } |w[f,f] | =
 \| \varphi\|_{ L^\infty ({\Bbb R}  )}.
\label{eq:BHb}\end{equation}
      \end{theorem}
     
 Estimate \e{eq:QFq1} means that there exists a bounded operator $W$ such that 
 $ w[f,g] =  ( f, Wg)$ for all $f, g\in L^2 ({\Bbb R}_{+} )$.   So Theorem~\ref{B-H} is   quite similar to the classical Toeplitz  result stating that the operator \e{eq:HDT} is bounded if and only if $t_{n}$ are the Fourier coefficients of a bounded function on the unit circle. Since we were not able to find a proof of      Theorem~\ref{B-H} in the literature,   it will be given  in Subsection~2.2 for completeness of our presentation.   However our main concern is to treat unbounded  Wiener-Hopf   operators.  
     
      The theory of Wiener-Hopf and Toeplitz  operators is a very well developed subject. 
 We refer to the books \cite{RR} (Chapter~3),   \cite{GGK} (Chapter~XII),  \cite{NK} (Chapters~B.4 and B.6),   \cite{Pe} (Chapter~3) and \cite{Bo}   for basic information on this theory.   However  results on unbounded Wiener-Hopf operators are practically nonexistent.

 \medskip

{\bf 1.2.}
In this paper, we consider semibounded Wiener-Hopf  operators $W$ in the space $L^2({\Bbb R}_{+})$ and follow basically the scheme of  \cite{YaT} where semibounded  Toeplitz operators were studied. However analytically this paper and \cite{YaT} are rather different.  

{\it We always suppose that $w(x)= \ov{w(-x)}$ so that the quadratic form \e{eq:QFq} is real and assume that
    \begin{equation}
w[f,f] \geq \gamma \| f\|^2  ,\q f\in C_{0}^\infty ({\Bbb R}_{+}), \q \| f\| =\| f\|_{L^2({\Bbb R}_{+})},
 \label{eq:T1}\end{equation}
 for some } $\gamma\in \Bbb R$.
In this   case, we are tempted to define $W$ as a self-adjoint operator corresponding to the quadratic form $w[f,f] $. Such an   operator exists if  the form $w[f,f] $ is closable in the space $L^2 ({\Bbb R}_{+})$, but as is well known this is not always true (for example, if $w(x)=1$ for all $x\in {\Bbb R}$). We refer to the book \cite{BSbook} for  basic information concerning these notions; they are also briefly discussed in Subsection~2.1.    We recall that, by  definition, the operator corresponding to the form $w[f,f]+ \beta \| f\|^2 $ is given by the equality $W_\beta = W+ \beta I$ (observe that the identity operator $I$ is a Wiener-Hopf operator).
Also by  definition, if a form $w[f,f]$  is closable, then all forms $w[f,f] +\beta \| f\|^2$  are closable. Therefore we can suppose that the number $\gamma$ in \e{eq:T1} is positive; for definiteness, we choose $\gamma=1$.
 
  We proceed from  the Bochner-Schwartz theorem.

  \begin{theorem}\label{R-H} 
  Let the form $w[f,f] $ be defined by the relation \e{eq:QFq} where the distribution $w\in C_{0}^\infty({\Bbb R})'$.  Then the condition
      \begin{equation}
w [f,f]  
 \geq 0, \q \forall f\in C_{0}^\infty ({\Bbb R}_{+}),
 \label{eq:QF}\end{equation}
  is satisfied if and only if there exists a
  non-negative   measure $d {\sf M} (\lambda)$ on the line $\Bbb R$ such that 
  \begin{equation}
w(x) = \frac{1}{2\pi} \int_{\Bbb R}  e^{-ix \lambda} d{\sf M}(\lambda ).
 \label{eq:WH}\end{equation}
Here the measure obeys the condition
    \begin{equation}
\int_{\Bbb R} (1+\lambda^2)^{-p} d{\sf M}( \lambda)<\infty
\label{eq:Sch}\end{equation}
for some $p $ $($that is,  it has at most a polynomial growth at infinity$)$.
    \end{theorem}

    Note  that usually    (see, e.g., Theorem~3 in \S 3 of Chapter II of the book \cite{GUEVI}) instead of condition \e{eq:QF} one requires that
    \begin{equation}
{\bf w} [f,f] : =\int_{\Bbb R} \int_{\Bbb R}   w(x-y) f(y) \ov{f (x)} dx dy 
 \geq 0, \q \forall f\in C_{0}^\infty ({\Bbb R}) ,
 \label{eq:QFx}\end{equation}
which   looks more restrictive. However the form ${\bf w} [f,f] $ is invariant with respect to shifts, that is, ${\bf w} [f_t,f_t] = {\bf w} [f ,f ] $ if $f_t (x)= f(x-t)$, and $ {\bf w} [f ,f ] =w [f,f]$ if $f\in C_{0}^\infty ({\Bbb R}_{+})$. Thus the conditions 
\e{eq:QF} and \e{eq:QFx} are equivalent.

Observe that the Lebesgue measure $d{\sf M}( \lambda)=d\lambda$ satisfies the condition \e{eq:Sch} with $p>1/2$. 
  So without  loss of generality, it is convenient to assume that
  $p>1/2$.   For  the Lebesgue measure,  relation \e{eq:WH} yields $w(x)=\d(x)$ (the delta-function) so that  $W=I$ and  $w[f,f] =  \| f\|^2$.  Therefore the measure corresponding to the form
    $w[f,f]+ \beta \| f\|^2 $ equals $d{\sf M}(\lambda ) + \beta d \lambda$, and  relation \e{eq:WH} extends to all semibounded Wiener-Hopf quadratic forms.   Thus we have the one-to-one correspondence between Wiener-Hopf  quadratic forms satisfying estimate  \e{eq:T1} and real measures satisfying the condition $ {\sf M}  (X)\geq \gamma    |X|$ ($  |X|$ is the Lebesgue measure) for all Borelian sets $X\subset{\Bbb R}$.

    We emphasize that a priori we only require that $w\in C_{0}^\infty({\Bbb R})'$, but, according to Theorem~\ref{R-H}, the semiboundedness condition \e{eq:T1} (and, in particular, \e{eq:QFq1}) ensures that $w\in {\cal S}'$.

  \medskip

{\bf 1.3.} 
  Our goal is to find necessary and sufficient conditions for the form $w[f,f]$ to be  closable. The answer to this question is strikingly  simple.
   
    \begin{theorem}\label{T1} 
    Let the form $w[f,f]$ be given by formula \e{eq:QFq} on elements $f\in C_{0}^\infty ({\Bbb R}_{+})$, and let the condition 
 \e{eq:T1} be satisfied  for some  $\gamma\in \Bbb R$.
   Then the form $w[f,f] $ is closable in the space $L^2 ({\Bbb R}_{+})$ if and only if  the measure $d{\sf M}  (\lambda)$ in the equation \e{eq:WH} is  absolutely continuous.
    \end{theorem}
    
 We always understand the absolute continuity with respect to the Lebesgue measure.   Therefore  Theorem~\ref{T1} means that $d{\sf M}  (\lambda)=  \varphi (\lambda) d\lambda$ where $\varphi \in L^1_{\rm loc} (\Bbb R)$,
   \begin{equation}
\int_{\Bbb R} (1+\lambda^2)^{-p} | \varphi (\lambda)|d  \lambda<\infty
 \label{eq:QFac}\end{equation}
       and $\varphi (\lambda) \geq \gamma$.  
       The function $\varphi(\lambda)$ is known as the symbol of the Wiener-Hopf operator $W$. Thus Theorem~\ref{T1} shows that in the semibounded case, the symbol of a  Wiener-Hopf operator can be correctly defined if and only if the corresponding quadratic form is closable.  So Theorem~\ref{T1} extends Theorem~\ref{B-H} from  bounded to semibounded operators.

 Our proof of Theorem~\ref{T1}   requires a
continuous analogue of the classical Riesz Brothers theorem. Let us state this result here. For a measure $d M ( \lambda)$ on $\Bbb R$, we denote by $d | M| ( \lambda)$ its  variation.

     \begin{theorem}\label{RBT} 
     Let   $d M(\lambda)$ be a complex     measure on the line $\Bbb R$ such that
      \begin{equation}
\int_{\Bbb R} (1+\lambda^2)^{-p} d | M| ( \lambda)<\infty
\label{eq:Schs}\end{equation} 
for some $p$. Put
  \begin{equation}
\sigma(x) = \frac{1}{2\pi} \int_{\Bbb R} e^{-ix \lambda} dM(\lambda ) 
 \label{eq:WHs}\end{equation}
 and suppose that $\sigma \in L^2 (a,\infty)$ for some $a\in \Bbb R$. Then the measure $d M(\lambda)$  is  absolutely continuous.
    \end{theorem}   
    
    We allow $a\in \Bbb R$ in   Theorem~\ref{RBT} to be arbitrary since, for example, the function $\sigma(x)=\d (x-x_{0})$ for any $x_{0}\in \Bbb R$ does not belong to $L^2_{\rm loc} ({\Bbb R})$, but the corresponding measure $d M(\lambda)= e^{i x_{0}\lambda} d\lambda$ is of course absolutely continuous.

    Theorems~\ref{T1}  and \ref{RBT}  extend naturally to vectorial Wiener-Hopf  operators $W$ and operator valued measures $d M(\lambda)$, but we do not dwell upon it here.

 \medskip

{\bf 1.4.}  
Section~2 is of a preliminary nature. In particular, we give here precise definitions of the quadratic forms \e{eq:QFq} and \e{eq:QFx}. The form ${\bf w}[f,f]$ is considered in Section~3 where we establish    a simplified version (Theorem~\ref{T1B}) of our main result,  Theorem~\ref{T1}. The substantial  difference between Theorems~\ref{T1} and \ref{T1B} is that, to consider Wiener-Hopf quadratic form, we need the continuous version of Riesz  Brothers theorem (Theorem~\ref{RBT}).  Theorems~\ref{RBT} and \ref{T1}  are proven in Sections~4 and 5, respectively. 

     A comparison of   our results   with  similar statements for Toeplitz operators are  postponed until Section~6.    There, we also discuss  a certain parallelism between theories of Wiener-Hopf and Hankel operators.

      \section{Wiener-Hopf operators. Generalities} 
    
     {\bf 2.1.}
     Let us  first briefly recall the notion of closable forms.
Let $w[f,f]$ be a quadratic form defined on a set $\cal D$ dense in a Hilbert space $\cal H$ and satisfying inequality \e{eq:T1} where $\| f\|$ is the norm of  $f\in\cal H$. Suppose that $\gamma=1$,  consider the norm $\| f\|_{W}= \sqrt{w[f,f]}$ and introduce the closure ${\cal D}[w]$ of $\cal D$  in this norm. If ${\cal D}[w]$ can be realized as a subset of $\cal H$, then one says that $w[f,f]$ is closable in the space $\cal H $; it means that the conditions
      \[
\| f_{n}\|\to 0 \q{\rm and}  \q \| f_{n}-f_{m}\|_W \to 0
\]
 as $n,m\to\infty$ imply that $\| f_{n} \|_{W}\to 0$. It is easy to show (see \S 10.3  of the book \cite{BSbook}) that if $W_{0}$ is a symmetric semibounded operator on $\cal D$, then the form $ w[f,f]: =  ( f, W_{0} f)$ is closable.
 
 Let the form $w[f,f]$ be closable. Then  it extends by continuity to all $f\in {\cal D}[w] \subset \cal H$, and one says that
 the form $w[f,f]$ is closed on   ${\cal D}[w] $.  For a closed form there exists a unique self-adjoint operator $W$ such that $W\geq I$ and
        \begin{align*}
        w[f,g] &=  ( f, Wg), \q \forall f\in {\cal D}[w], \q \forall g\in  {\cal D} (W)\subset {\cal D}[w],
        \\
        w[f,f]&= \| \sqrt{W} f\|^2, \q \forall f\in  {\cal D}(\sqrt{W}) = {\cal D}[w].
            \end{align*}
    Note that the domain $ {\cal D} (W) $ of the operator $W$  does not admit an efficient description.

 We are  going to use these general definitions for the space ${\cal H}  = L^2 ({\Bbb R}_{+})$ and the Wiener-Hopf  quadratic forms \e{eq:QFq} on $\cal D=C_{0}^\infty ({\Bbb R}_{+})$ or  for the space ${\cal H}  = L^2 ({\Bbb R})$ and the  convolution quadratic forms \e{eq:QFx} on $\cal D=C_{0}^\infty ({\Bbb R})$.
 
 Of course  quadratic forms, in particular, the  Wiener-Hopf  forms, are not necessarily  closable.
 
     \begin{example}\label{ex}
Let $w(x)=1$ for all $x\in {\Bbb R}$. Adding the term $\| f\|^2$,  we obtain the form
    \[
w[f,f]= \big|  \int_{{\Bbb R}_{+}}  f(x) dx\big|^2 +    \int_{{\Bbb R}_{+}}  \big|  f(x) \big|^2 dx  
 \]
satisfying inequality \e{eq:T1}  with $\gamma=1$.
Let $\psi\in C_{0}^\infty ({\Bbb R}_{+})$ and $\int_{0}^\infty \psi(x)dx=1$.
 Define the sequence $f_{n}\in C_{0}^\infty ({\Bbb R}_{+})$ by the equalities $f_{n}(x)=n^{-1}\psi(n^{-1}x)$. Then  $\| f_{n}\|=n^{-1/2}\to 0$ as $n\to\infty$. Since $\int_{0}^\infty f_{n}(x)dx=1$ for all $n$, we  have $ \| f_{n}- f_{m}\|_{W}= \| f_{n}- f_{m}\|\to 0$ as $n,m\to\infty$.  Nevertheless  $ \| f_{n} \|_{W} \geq 1$.
    \end{example}
    
    Note that the   measure $d {\sf M}(\lambda)$ corresponding in  \e{eq:WH} to the function  $w(x)=1$, $\forall x\in {\Bbb R}$, is supported by the point $0  $: ${\sf M} (\{0\})=2\pi $, ${\sf M} ({\Bbb R} \setminus \{0\})=0$.
    
    On the other hand, we have  the following simple assertion.
    
       \begin{lemma}\label{ex1}
If    $w \in L^2 ({\Bbb R})$ and  the form \e{eq:QFq} is semibounded, then it is closable.
    \end{lemma}
    
\begin{pf}
Under the assumption $w \in L^2 ({\Bbb R})$,  the    operator \e{eq:HD} (it will be denoted
$W_{0}$) is correctly defined   on  $f\in C_{0}^\infty ({\Bbb R}_{+})$. This operator is
 symmetric  because $w(x)=\ov{w(-x)}$.
Since  $w[f,f]=  ( f,  W_{0} f)$,   the form $w[f,f]$ is closable.
    \end{pf}
    
    Note that under the assumptions of Lemma~\ref{ex1} the self-adjoint operator $W$ corresponding to the form $w[f,f]$ is the Friedrichs extension of $W_{0}$.
    
\medskip

    {\bf 2.2.}    Let us now discuss the precise definitions of the Wiener-Hopf and convolution quadratic forms \e{eq:QFq}
    and \e{eq:QFx}.  Obviously, \e{eq:QFx} can be written as
     \begin{equation}
{\bf w}[f,f] = \la w, f\circ \bar{f}   \ra
 \label{eq:QFq3}\end{equation}
 where 
   \begin{equation}
(f\circ g) (x)= 
\int_{\Bbb R}     f(y - x)   g (y) dy 
 \label{eq:QFq2}\end{equation}
 is the convolution composed with the reflection of $f$ and $g$;     the duality symbol $\la \cdot,  \cdot \ra$ is induced by the complex  scalar product in $L^2 ({\Bbb R})$:
 \[
 \la w,  \theta \ra=\int_{\Bbb R} w(x) \ov{\theta (x)} dx
 \]
 (so $w$ is an antilinear functional).
   Since $f\circ g\in C_{0}^\infty ({\Bbb R})$ for $f , g\in C_{0}^\infty ({\Bbb R})$, the form 
 \e{eq:QFq3} is well defined if $w\in C_{0}^\infty ({\Bbb R})'$ for all $f  \in C_{0}^\infty ({\Bbb R})$ and, in  particular, for
 $f  \in C_{0}^\infty ({\Bbb R}_{+})$. We now put
   \begin{equation}
  w [f,f] := {\bf w}[f,f] \q {\rm for}\q f  \in C_{0}^\infty ({\Bbb R}_{+}).
 \label{eq:WHC}\end{equation}
  We usually write $w[f,f] $ and  ${\bf w}[f,f]$ in the forms \e{eq:QFq} and \e{eq:QFx} keeping in mind that their precise definitions are given by formulas \e{eq:QFq3}, \e{eq:QFq2} and \e{eq:WHC}.

 For $t \in{\Bbb R}$, we put $f_t (x)= f(x- t)$. Since $f_t\circ g_t =f\circ g$, we see that the form ${\bf w} [f,f] $ is invariant with respect to shifts, that is,    ${\bf w} [f_t,f_t ] = {\bf w} [f ,f ] $. Therefore relation \e{eq:WHC} implies that the conditions 
\e{eq:QFq1} and 
 \begin{equation}
|{\bf w}[f,f] | \leq C \| f\|^2 ,   \q \| f\|=\| f\|_{L^2 ({\Bbb R})} , \q \forall f \in C_{0}^\infty ({\Bbb R}),
 \label{eq:QFq1x}\end{equation}
  are equivalent.
  
  We standardly define the Fourier transform 
     \[
(\Phi u) (x)= \frac{1}{\sqrt{2\pi }}
\int_{\Bbb R}   e^{-ix\lambda}  u(\lambda) d\lambda.
 \]
 Of course the operator $\Phi : L^2 ({\Bbb R}) \to L^2 ({\Bbb R})$ is unitary and $\Phi : {\cal S} \to {\cal S}$, $\Phi : {\cal S}' \to {\cal S}'$. Below $C$ are positive constants whose precise values are of no importance.
 
 Now we are in  a position to prove Theorem~\ref{B-H}. 
 
   \begin{pf}
   As  already explained, estimate \e{eq:QFq1} for $f\in C_{0}^\infty ({\Bbb R}_{+})$ implies estimate \e{eq:QFq1x} for
   $f\in C_{0}^\infty ({\Bbb R})$. Moreover, we can extend it to sesquilinear forms which yields
     \begin{equation}
|\la w, f\circ \bar{g}   \ra | =| w [f,g] |    \leq C \|f\|  \| g\|, \q \forall f,g\in C_{0}^\infty := C_{0}^\infty ({\Bbb R}).
 \label{eq:BF}\end{equation}
 We use that
 \begin{equation}
(\Phi^{*} (f\circ \bar{g}) ) (\lambda)= \sqrt{2\pi } (\Phi^{*} f )(-\lambda)  (\Phi^{*} \bar{g}) (\lambda)
 \label{eq:BF1}\end{equation}
 where both factors on the right belong to the space $\wh{C_{0}^\infty}: =\Phi^{*} C_{0}^\infty$.
 Let $\wh{w}= \Phi^{*} w\in \wh{C_{0}^\infty}'$, $f_{1} (x)= f(-x)$, $g_{1} (x)= \ov{g(x)}$. Then it follows from \e{eq:BF} and \e{eq:BF1} that
   \begin{equation}
\sqrt{2\pi } |\la \wh{w}, \wh{f}_{1} \,\wh{g}_{1}   \ra | = |\la w, f\circ \bar{g}   \ra | \leq C \|\wh{f}_{1}\|  \| \wh{g}_{1}\| 
 \label{eq:BF2}\end{equation}
 forall $\wh{f}_{1}=\Phi^* f_{1}\in \wh{C_{0}^\infty }$ and $\wh{g}_{1}=
\Phi^* g_{1}\in \wh{C_{0}^\infty }$.
 
 For every $u\in L^{1} ({\Bbb R})$, we set     $F(\lambda)=\sqrt{|u(\lambda)|}$ and   $G(\lambda)=   u (\lambda) F(\lambda)^{-1}$. Then   $u=FG$  and
 \[
   \|F\|_{L^{2}({\Bbb R})}^2 = \|G\|_{L^{2}({\Bbb R})}^2 =  \|u\|_{L^{1}({\Bbb R})}.
 \]
 Approximating in $L^{2}({\Bbb R})$ the functions $F$ and $G$ by functions in 
  $\wh{C_{0}^\infty}$, we see that estimate \e{eq:BF2} implies the  estimate 
  \begin{equation}
|\la \wh{w}, u   \ra | \leq C    \|u\|_{L^{1}({\Bbb R})} 
 \label{eq:BF3}\end{equation}
 on a set of functions $u$ dense in the space $L^{1}({\Bbb R})$. Its dual is $L^{\infty}({\Bbb R})$, and hence it follows from \e{eq:BF3} that
   \[
\la \wh{w}, u   \ra = (2\pi)^{-1/2}  \int_{\Bbb R} \varphi (\lambda) \ov{u(\lambda)} d\lambda
\]
 for some $\varphi\in L^{\infty}({\Bbb R})$. Therefore $\wh{w}=(2\pi)^{-1/2}  \varphi$, which yields  formula  \e{eq:BH}.
  
    Conversely, if \e{eq:BH}  is satisfied, then 
       \[
\la w, f\circ \bar{f}   \ra  =  
\int_{\Bbb R}  \varphi (\lambda) | \wh{f}(\lambda)|^2 d\lambda
 \]
at least for all  functions $f\in C_{0}^\infty$.
This implies equality  \e{eq:BHb}.
        \end{pf}

      \begin{remark}\label{www}
 As was already explained, estimates \e{eq:QFq1} and \e{eq:QFq1x} are equivalent. So under the assumptions of Theorem~\ref{B-H} the   operator ${\bf W}$ corresponding to the quadratic form ${\bf w}[f,f]$ is bounded.
     \end{remark}
    
The convolution operator ${\bf W}$ and the Wiener-Hopf operator $W$ satisfy the relation
  \begin{equation}
({\bf W} f,f)= (W   f,   f), \q \forall f \in L^{2} ({\Bbb R}_{+} ),
 \label{eq:WW}\end{equation}
 (of course $f(x)$ is extended by zero to $x<0$).  This relation defines the Wiener-Hopf operator $W$ in terms of the convolution operator ${\bf W}$.
 
\medskip

    {\bf 2.3.} 
     As in the discrete case, the Wiener-Hopf operators in  the space $L^2 ({\Bbb R}_{+})$ can be characterized by a commutation relation. Let $S_t$, $t\geq 0$,  be the shift in the space $L^2 ({\Bbb R}_{+})$:
  \[
(S_t f) (x)= f(x- t) \q {\rm for}  \q x\geq  t  \q {\rm and} \q (S_t f) (x)= 0 \q {\rm for} \q x < t. 
\]
 Obviously, $S_t f\in C_{0}^\infty ({\Bbb R}_{+})$ if $  f\in C_{0}^\infty ({\Bbb R}_{+})$ and, by definition  \e{eq:QFq2}, we have $S_t  f\circ S_t g =   f\circ  g$ for all $  f,g \in C_{0}^\infty ({\Bbb R}_{+})$  and $ t \geq 0$.  Thus it follows from   \e{eq:QFq3}   that
     \[
w[S_t f, S_tf] = w[f,f] , \q \forall t\geq 0, 
\]
 so that bounded Wiener-Hopf operators   satisfy the commutation relation
  \begin{equation}
S_t ^* W S_t = W, \q \forall t\geq 0.
 \label{eq:FS1}\end{equation}
 
 The converse statement is also true. It is not a direct consequence of the corresponding result for Toeplitz operators, but its proof   essentially follows the same scheme (see, e.g., the book \cite{Pe}, Chapter~3, Theorems~1.1 and 1.2).
 
   \begin{theorem}\label{CR}
  If a bounded operator $W$ in the space $L^{2} ({\Bbb R}_{+})$ satisfies relation \e{eq:FS1}, then $W$ is a Wiener-Hopf operator, that is, its quadratic form is  given by  the equality $(Wf,f)= \la w, f\circ\bar{f}\ra$  where the distribution   $w $ admits representation  \e{eq:BH} with $\varphi\in L^\infty ({\Bbb R})$.
     \end{theorem}
 
   \begin{pf}
   Let ${\bf S}_{a}$ be the shift in the space $L^2 ({\Bbb R})$ defined by the formula
   \[
   ({\bf S}_{a} f) (x)= f(x-a), \q a\in {\Bbb R}.
   \]
   Let ${\cal H}_{T}\subset L^2 ({\Bbb R})$ consist of functions $f(x)$ such that $f(x)=0$ for $x\leq -T$. Clearly, the union $\sf D$ of the sets ${\cal H}_{T}$ over all $T>0$ is dense in $L^2 ({\Bbb R})$. Given  a bounded Wiener-Hopf operator $W$ in $L^{2} ({\Bbb R}_{+})$,   we  introduce the convolution operator $\bf W$ in  $L^2 ({\Bbb R})$ by its quadratic form defined on $\sf D$: 
    \begin{equation}
 ({\bf W} f,f)= (  W {\bf S}_{a} f, {\bf S}_{a}f), \q  a\geq T \q {\rm if}\q f\in {\cal H}_{T}.
 \label{eq:WW1}\end{equation}
 Obviously, $ {\bf S}_{a}f\in L^2 ({\Bbb R}_{+})$ and in view of commutation relation \e{eq:FS1}
  \[
  (  W {\bf S}_{a} f, {\bf S}_{a}f)=   (  W {\bf S}_{a-T} {\bf S}_{T}  f, {\bf S}_{a-T} {\bf S}_{T} f)=   (  W   {\bf S}_{T}  f,   {\bf S}_{T} f)
\]
 so that the right-hand side of \e{eq:WW1} does not depend on $a\geq T $. It  follows from \e{eq:WW1} that $\| {\bf W}\|\leq \| W\|$.

 Let us now check that
  \begin{equation}
 ({\bf W} f,f)= (  {\bf W} {\bf S}_{b} f, {\bf S}_{b}f), \q \forall b\in {\Bbb R}.
 \label{eq:WW3}\end{equation}
   We may suppose that $f\in {\cal H}_{T}$ for some $T>0$. Put  $a = \max\{T, T-b\}$. By definition \e{eq:WW1}, relations \e{eq:WW3} and
 \begin{equation}
(  W {\bf S}_{a} f, {\bf S}_{a}f)= (  W {\bf S}_{b} {\bf S}_{a}f, {\bf S}_{b}{\bf S}_{a}f) 
 \label{eq:WW4}\end{equation}   
      are equivalent. If $b>0$, then   ${\bf S}_{a} f ={\bf S}_{T} f \in {\cal H}_0\cong L^2 ({\Bbb R}_{+})$,  so that equality \e{eq:WW4} follows from \e{eq:FS1} for $t=b$. If $b<0$, then $a= T-b$ and \e{eq:WW4} can be rewritten as
       \[
(  W {\bf S}_{-b} {\bf S}_{T} f, {\bf S}_{-b} {\bf S}_{T}f)= (  W {\bf S}_{T} f, {\bf S}_{T} f) .
\] 
 Since ${\bf S}_{T} f \in {\cal H}_0 $, this equality follows again from \e{eq:FS1} for $t=-b$.
 
 Let us now set $\wh{{\bf W} }= \Phi^{*}{\bf W} \Phi$ and $\wh{{{\bf S} } }_{b}= \Phi^{*}{\bf S}_{b} \Phi$. According to \e{eq:WW3} we have
     \begin{equation}
\wh{{\bf W} }\wh{{{\bf S} } }_{b}= \wh{{{\bf S} } }_{b}\wh{{\bf W} }, \q \forall b\in {\Bbb R} .
 \label{eq:WW6}\end{equation}
 Since the operator $\wh{{{\bf S} } }_{b}$ acts as multiplication by $e^{ib\lambda}$, relation \e{eq:WW6} implies that 
  the operator $\wh{{\bf W} }$ acts as multiplication by a function $\varphi(\lambda)$. The function $\varphi\in L^\infty ({\Bbb R})$ because the operator $\wh{{\bf W} }$ is bounded. Therefore ${\bf W}$ is the convolution operator with integral  kernel $w(x)$ satisfying \e{eq:BH}. Since  relation \e{eq:WW1} with $a=T=0$ yields \e{eq:WW}, we obtain all the conclusions about the  operator $W$.
       \end{pf}

  \medskip
  
  \section{Semibounded convolution operators}

In this section we prove a simplified version of Theorem~\ref{T1} where the Wiener-Hopf quadratic form  \e{eq:QFq} is replaced by the convolution quadratic form  ${\bf w} [f,f]$.

\medskip

{\bf 3.1.}
Let the   quadratic form ${\bf w} [f,f]$ be defined by formula  \e{eq:QFq3} where $w\in  C_{0}^\infty ({\Bbb R})'$ and $f\in C_{0}^\infty ({\Bbb R})$. As before we  suppose that $w(x)= \ov{w(-x)}$ so that this quadratic form  is real and assume that
    \begin{equation}
{\bf w} [f,f] \geq \gamma \| f\|^2  ,\q f\in C_{0}^\infty ({\Bbb R}), \q \| f\| =\| f\|_{L^2({\Bbb R})},
 \label{eq:T1B}\end{equation}
 for some $\gamma\in \Bbb R$. Then the representation \e{eq:WH}  is satisfied with the measure $d{\sf M} (\lambda)$ obeying condition \e{eq:Sch} and such that ${\sf M} (X)\geq \gamma |X|$  for all $X\subset {\Bbb R}$.

   Our goal is to find necessary and sufficient conditions for the form ${\bf w}  [f,f]$ to be  closable.

 \begin{theorem}\label{T1B} 
    Let the form ${\bf w}[f,f]$ be given by formula \e{eq:QFq3} on elements $f\in C_{0}^\infty ({\Bbb R} )$, and let the condition 
 \e{eq:T1B} be satisfied. 
   Then the form ${\bf w} [f,f] $ is closable in the space $L^2 ({\Bbb R} )$ if and only if  the measure $d{\sf M}  (\lambda)$ in the equations \e{eq:WH} is  absolutely continuous. In this case, it is closed on all functions $f\in L^2 ({\Bbb R} )$ such that
   \begin{equation}
   \int_{\Bbb R}  | \wh{f} (\lambda )|^2 d {\sf M}(\lambda)<\infty,\q \wh{f}=\Phi^{*} f.
 \label{eq:Ab}\end{equation}
    \end{theorem}
         
  \medskip
  
  {\bf 3.2.}  
  By the proof of Theorem~\ref{T1}  we may suppose that estimate    \e{eq:T1B}  is true for $\gamma=1$. Then  the equations  \e{eq:WH} are satisfied with a measure $d{\sf M}  (\lambda )$ such  that ${\sf M} (X)\geq |X|$ for all Borelian sets  $X\subset \Bbb R$;  in particular, the measure $d {\sf M}  (\lambda)$  is positive. 
  
        Our proof relies on the following auxiliary construction. Let $   L^2 (\Bbb R; d{\sf M} )$ be the space of functions $u (\lambda)$ on $\Bbb R$ with the norm 
   \[
   \| u\|_{ L^2 (\Bbb R; d {\sf M} )}=\sqrt{\int_{\Bbb R} | u(\lambda)|^2 d{\sf M} (\lambda)}.
   \]
       We define  an operator $A \colon L^2 ({\Bbb R})\to    L^2 (\Bbb R; d {\sf M} )$  on domain $\cal D (A)=  C_{0}^\infty ({\Bbb R})$ by the formula
   \begin{equation}
(Af) (\lambda)=\frac{1} {\sqrt{2\pi}}\int_{\Bbb R} e^{ix\lambda} f(x) dx
 \label{eq:A}\end{equation}
 Since         ${\cal S} \subset L^2 (\Bbb R; d {\sf M} )$ according to \e{eq:Sch}, we see that
    $ A f\in {\cal S} \subset L^2 (\Bbb R; d {\sf M} )$ for all $f \in C_{0}^\infty ({\Bbb R})$.
 Obviously,  the operator $A$ acts the (inverse) Fourier transform, but it is considered as a mapping of $ L^2 ({\Bbb R})$ into $L^2 (\Bbb R; d {\sf M} )$.
 In view of  equation  \e{eq:WH} the form       \e{eq:QFx}   can be written as
   \begin{multline}
{\bf w} [f,f] =\frac{1} { 2\pi }\int_{\Bbb R}\int_{\Bbb R} dx dy f(y) \ov{f(x)} \Big(\int_{\Bbb R}e^{-i (x-y)\lambda}  d {\sf M}(\lambda) \Big)
\\
    = \frac{1} {2\pi} \int_{\Bbb R}  \big|  \int_{\Bbb R}e^{i x\lambda} f(x)dx\big|^2 d {\sf M}(\lambda) =  \| Af\|^2_{ L^2 (\Bbb R; d{\sf M})} , \q f\in C_{0}^\infty ({\Bbb R}).
 \label{eq:A2}\end{multline}
 We have interchanged the order of integrations in $x,y$ and $\lambda$ here. Of course the Fubini theorem is not applicable  now. Nevertheless equality \e{eq:A2} is true because $f\in {\cal S}$ and the Fourier transforms are understood in the sense of the Schwartz space ${\cal S}'$.

 Equality \e{eq:A2} yields the following result.
 
    \begin{lemma}\label{de}
 The form $ {\bf w} [f,f]$ defined   on $C_{0}^\infty ({\Bbb R})$   is closable in the space $L^2 ({\Bbb R} )$  if and only if the operator $A \colon L^2 ({\Bbb R} )\to    L^2 (\Bbb R; d{\sf M})$ defined   on the  domain $\cal D (A)=C_{0}^\infty ({\Bbb R})$  by formula \e{eq:A} is closable.  In this case the form $ {\bf w} [f,f]$ is closed on the domain $\cal D [ {\bf w} ]=\cal D (\cl A)$ of the closure   of the operator $A$.
    \end{lemma}

  Now it is easy to prove the `` if " part of Theorem~\ref{T1B}.
  Suppose   that the measure $d{\sf M}  (\lambda)$ is absolutely continuous, that is, 
      \begin{equation}
  d{\sf M}  (\lambda)= \varphi (\lambda) d\lambda
   \label{eq:AC}\end{equation}
   where $ \varphi (\lambda) \geq 1$ and, for some $p$, condition \e{eq:QFac} is satisfied. Let us check that then  the form $ {\bf w} [f,f]$ defined   on $C_{0}^\infty ({\Bbb R})$ (or on the Schwartz class $\cal S$)  is closable in the space $L^2 ({\Bbb R} )$. In view of Lemma~\ref{de} it suffices to  verify the same fact for the operator $A $.
   It follows from \e{eq:A} and \e{eq:AC}  that
       \begin{equation}
\| Af\|^2_{ L^2 (\Bbb R; d{\sf M})}
     = \int_{\Bbb R}  | \wh{f}(\lambda) |^2 \varphi (\lambda)d \lambda  , \q \wh{f}= A f.
      \label{eq:A2z}\end{equation}
 Thus \e{eq:A2z}  is the quadratic form of the operator of multiplication by $ \varphi (\lambda)$ in the space $L^2 ({\Bbb R})$ defined on     $\cal S$.
It is closable because    $\varphi\in L^1_{\rm loc} ({\Bbb R}  )$. Moreover,       the form \e{eq:A2z}   is closed on the set of all  $\wh{f}\in L^2({\Bbb R} )$ such that the integral \e{eq:A2z} is finite. So the operator $A$ defined on $C_{0}^\infty ({\Bbb R})$  (or on   $\cal S$)  is closable, and $f\in {\cal D} (\cl A)$ if and only if integral \e{eq:A2z} is finite. Lemma~\ref{de} allows us to carry over these results to the form ${\bf w} [f,f]$.

     \medskip
  
  {\bf 3.3.} 
  To prove the converse statement, we have
    to construct the  operator $A^* \colon L^2 ({\Bbb R; d{\sf M}} )\to    L^2 (\Bbb R )$  adjoint to the operator $A \colon L^2 ({\Bbb R} )\to    L^2 (\Bbb R; d{\sf M})$.  Since    
 \[
\big| \int_{\Bbb R} \ov { f(\lambda)} u(\lambda) d {\sf M} (\lambda)\big|^2\leq
\int_{\Bbb R} | f (\lambda)|^2   d {\sf M} (\lambda) \int_{\Bbb R} |u (\lambda)|^2   d {\sf M} (\lambda)<\infty 
\]
 for an arbitrary $u\in L^2 (\Bbb R; d {\sf M})$ and all $f \in {\cal S}$,
 the distribution $u(\lambda) d {\sf M} (\lambda)$ belongs to the class $ {\cal S}'$. Therefore its Fourier transform
    \begin{equation}
 u_{*} (x) = \frac{1} {\sqrt{2\pi}} \int_{\Bbb R} e^{-ix\lambda} u(\lambda) d {\sf M} (\lambda)
 \label{eq:A1}\end{equation}
is correctly defined (as usual,   in the sense of the Schwartz space ${\cal S}'$) and $u_{*}\in{\cal S}'$.

     \begin{definition}\label{def}
 The set   ${\cal D}_{*}\subset   L^2 (\Bbb R; d {\sf M})$ consists  of all  $u \in  L^2 (\Bbb R; d {\sf M})$ such that   $u_{*} \in  L^2 ({\Bbb R} )$.
    \end{definition}

   \begin{lemma}\label{LTM}
 The operator $A ^*$ is given by the equality 
      \[
        (A^ {*}u) (x)=  u_{*} (x)
\]
  on the domain ${\cal D}(A^*) ={\cal D}_{*}$. 
    \end{lemma}
    
     \begin{pf}
     Obviously, for all $f\in \cal S$ and all $u\in L^2 (\Bbb R; d {\sf M})$, we  have the equality
       \begin{equation}
 (Af,u)_{L^2 (\Bbb R; d \sf  M)} =   \frac{1} {\sqrt{2\pi}} \int_{\Bbb R} \big( \int_{\Bbb R} e^{ix\lambda} f(x) dx\big) \ov{u(\lambda)}d {\sf M}(\lambda)=\int_{\Bbb R} f(x)   \ov{u_{*} (x)}    dx .
\label{eq:Y1}\end{equation} 
As usual, the Fourier transforms are here understood in the sense of ${\cal S}'$ so that the integrations over $x$ and $\lambda$ can be automatically interchanged. In particular, if 
$u\in {\cal D}_{*} $, then $u_{*} \in  L^2 ({\Bbb R} )$ and the right-hand side of \e{eq:Y1} equals $(f,u_{*})_{ L^2 ({\Bbb R} )}$. 
It follows that $   {\cal D}_{*}\subset {\cal D}(A^*)  $.

         Conversely, if $u \in {\cal D}(A^*) $, then 
         \[
         |(Af, u)_{L^2 (\Bbb R; d{\sf M})} | =   |(f, A^* u)_{ L^2 ({\Bbb R})} | \leq   \| A^* u\|_{ L^2 ({\Bbb R} )}  \,  \| f\|_{ L^2 ({\Bbb R} )}
             \]
          for all $ f \in {\cal S}$. Therefore it follows from equality \e{eq:Y1} that
               \[
\big|   \int_{\Bbb R} f(x)   \ov{u_{*} (x)}    dx \big| \leq  \| A^* u\|_{ L^2 ({\Bbb R} )}  \| f \|_{ L^2 ({\Bbb R} )}, \q \forall f \in {\cal S}.
\]
Since $ {\cal S}$ is dense in $ L^2 ({\Bbb R} )$, we  see that
 $ u_{*} \in L^2 ({\Bbb R} )$, and hence $u\in {\cal D}_{*}$.
          Thus   ${\cal D}(A^*)   \subset {\cal D}_{*}$.
          \end{pf}

      Recall that    an operator $A $   is closable if and only if its adjoint operator $A^*   $ is densely defined. Below $ \cl{\cal D}_{*}$ is the closure of the set ${\cal D}_{*}$
      in the space $L^2 ({\Bbb R}; d{\sf M}) $. So Lemmas~\ref{de} and \ref{LTM}  imply the following intermediary result.
      
       \begin{lemma}\label{adj}
 The operator $A  $ and the form ${\bf w}[f,f]$ are closable if and only if  
   \begin{equation}
 \cl{\cal D}_{*} =L^2 ({\Bbb R}; d{\sf M}). 
 \label{eq:D}\end{equation}
    \end{lemma}

    Now we are in a position to show that if   the form ${\bf w}[f,f]$ is closable, then the measure $d{\sf M}  (\lambda)$ is absolutely continuous. If $u\in {\cal D}_{*}$, then the function   \e{eq:A1} belongs to $L^2 ({\Bbb R})$ and hence, by  the Parseval theorem,
    \begin{equation}
 u(\lambda) d {\sf M} (\lambda)= \psi(\lambda) d\lambda
 \label{eq:Par}\end{equation}
 for some $\psi\in L^2 ({\Bbb R})$. In particular, the measure $ u(\lambda) d {\sf M} (\lambda)$ is absolutely continuous. 
    
 Next we use    the following simple assertion.  Below $\1_{X}$ is the characteristic function of a Borelian set $X\subset \Bbb R$.
    
      \begin{lemma}\label{ac}
      Suppose that a set ${\cal D}_{*}$ satisfies condition  \e{eq:D}. Let the measures $ u(\lambda) d {\sf M} (\lambda)$ be absolutely continuous
      for all $u\in{\cal D}_{*}$.  Then the measure $ d {\sf M} (\lambda)$  is also  absolutely continuous.
    \end{lemma}

  \begin{pf}  
  Let $|X|=0$.
    Suppose first that a set $X$ is bounded and hence $\1_{X}\in L^2 ({\Bbb R}; d{\sf M} )$. It follows from \e{eq:D} that there exists a sequence $u_{n}\in{\cal D}_{*}$ such that
  \[
\lim_{n\to\infty}  \| u_{n} -\1_{X} \|_{L^2 ({\Bbb R}; d{\sf M} )}=0,
\]
whence  
 \[
\big|  \int_{X} u_{n}(\lambda) d{\sf M} (\lambda) - {\sf M} (X) \big| \leq  \| u_{n} -\1_{X} \|_{L^2 ({\Bbb R}; d{\sf M} )} \sqrt{{\sf M} (X)} \to 0
\]
as $n\to\infty$.  Since the measures $ u_{n}(\lambda) d {\sf M} (\lambda)$  are absolutely continuous,  the integrals in the left-hand side are zeros so that ${\sf M} (X)=0$.
If a set $X$ is unbounded, then ${\sf M} (X\cap (-r,r)) =0$ for all $r<\infty$, and hence ${\sf M} (X ) =0$.
 \end{pf} 
 
 Now it is easy to conclude the ``only if " part of  Theorem~\ref{T1B}. 
  Suppose that the form ${\bf w} [f,f]$  is closable. Then by Lemma~\ref{adj} the condition  \e{eq:D} is satisfied. 
  By the definition of the set ${\cal D}_{*}$,    the measures \e{eq:Par} are absolutely continuous.      Hence   by Lemma~\ref{ac} the same is true for the measure $  d {\sf M} (\lambda)$. 
  
  It remains to show that  ${\cal D} ( \cl{A})$  consists of all $f\in L^2 ({\Bbb R} )$ such that condition  \e{eq:Ab} is satisfied. By definition,   $f\in {\cal D} ( \cl{A})$ if and only if there exists a sequence $f_{n}  \in C_{0}^\infty ({\Bbb R}) $ such that $f_{n}\to f$ in $L^2 ({\Bbb R})$ and $\wh{f}_{n}\to u= (\cl{A}) f$ in $L^2 ({\Bbb R}; d{\sf M})$. It follows that $\wh{f}_{n}\to \wh{f}$   and $\wh{f}_{n}\to u $ in $L^2 ({\Bbb R} )$ so  that  the function $\wh{f}=u\in L^2 ({\Bbb R}; d{\sf M})$. To prove the converse statement, observe that the kernel $\Ker (A^*)=\{0\}$. Indeed,  if integral \e{eq:A1}  is zero for a.e. $x\in {\Bbb R}$, then in view of \e{eq:Par} we have $u(\lambda) \psi(\lambda)=0$ and hence $u(\lambda)  =0$ for a.e. $\lambda\in {\Bbb R}$. Therefore for the image of $A$ we have
$  \cl \big(\Im(A)\big) =L^2 ({\Bbb R}; d{\sf M})$. Thus if condition  \e{eq:Ab} is satisfied for some $f\in L^2 ({\Bbb R} )$, then  there exists  a  sequence $f_{n}  \in C_{0}^\infty ({\Bbb R}) $ such that $A f_{n}= \wh{f}_{n}\to (\cl A)f= \wh{f}$ in $L^2 ({\Bbb R}; d{\sf M})$ and hence $f_{n}\to f$ in $L^2 ({\Bbb R})$. This means that $f\in{\cal D} ( \cl{A})$.  $\Box$

  \begin{remark}\label{Clo1} 
     Recall that the operator $A$ was defined by formula \e{eq:A}  on the domain ${\cal D} (A)= C_{0}^\infty ({\Bbb R}) $.
     Of course   $\cl{A}=A^{**}$ if $A$ is closable.  Let $A_{\rm max}$ be given  by the same formula \e{eq:A}    on the domain $\cal D(A_{\rm max})$ that consists of all $f\in L^2 ({\Bbb R} )$ such that $A_{\rm max} f\in L^2 (\Bbb R; d {\sf M})$. Then  the assertion of Theorem~\ref{T1B}  is equivalent to the equality
     $    \cl{A}= A_{\rm max}$.
        \end{remark}
        
 \section{Continuous analogue of the Riesz Brothers theorem} 
 
 Our goal in this section is to prove Theorem~\ref{RBT}.  
 
 \medskip
  
  {\bf 4.1.}    Let us proceed from the classical       Riesz Brothers theorem  (see, e.g., \cite{Hof}, Chapter~4).

         \begin{theorem}\label{brothers}
      Let $d\mu (z)$ be a  complex $($finite$)$ measure      on the unit circle $\Bbb T$. Suppose that its Fourier coefficients
              \[
 \int_{\Bbb T}  z^{-n}d\mu(z) =0\q {\rm for} \q n=1,2, \ldots.
\]
     Then  the measure $d\mu   (z)$ is absolutely continuous.
    \end{theorem}
    
    We  need this result in the following form. Let us introduce the measure $d\mu_{0} (z)$ on $\Bbb T$ supported by the point $1 $: $\mu_{0} (\{1\})=1$, $\mu_{0} (Y )=0$ if $Y \subset \Bbb T \setminus \{1\}$.
    
     \begin{corollary}\label{brothers1}
  Suppose that for some $b\in {\Bbb C}$ and some integer $N\geq 0$ the Fourier coefficients
    \begin{equation}
 \int_{\Bbb T}  z^{-n}d\mu(z) = b \q {\rm for} \q n=N+1, N+2, \ldots.
\label{eq:br}\end{equation}
     Then  the measure $d\mu   (z)- b d\mu_{0}   (z)$ is absolutely continuous. If, additionally, $\mu  (\{1\})=0$, then  $b=0$ and the measure $d\mu   (z) $ is also absolutely continuous.
    \end{corollary}
    
   \begin{pf}
   Set 
    \[
 d\mu_{1}(z)= z^{-N}  \big(d\mu (z)- b d\mu_0(z) \big).
\]
It follows from  \e{eq:br} that the measure $ d\mu_{1}(z)$ satisfies the condition of Theorem~\ref{brothers}, and hence it is absolutely continuous. The same is also true for the measure $d\mu   (z)- b d\mu_{0}   (z)$. If  $\mu  (\{1\})=0$, then $ b \mu_{0}   (\{1\})=0$ so that $b =0$.
      \end{pf}
      
   As we will see,   Theorems~\ref{RBT} and  \ref{brothers} are essentially equivalent if $p=0$ in condition \e{eq:Schs}, but Theorem~\ref{RBT} is more general if $p>0$ when $| M | (\Bbb R)=\infty$. However even in the general case,  we use Theorem~\ref{brothers} for the proof of Theorem~\ref{RBT}.
   
  Let us introduce an auxiliary measure
      \begin{equation}
 dm(\lambda )  =
   (1+\lambda^2)^{-p} dM(\lambda ).
 \label{eq:WHs2}\end{equation}
  According to condition \e{eq:Schs} it is finite, that is, $|m|({\Bbb R})<\infty$. Then we
   reduce     the problem on the line $\Bbb R$ to the problem on the circle $\Bbb T$ using the standard mapping $\omega :{\Bbb R}\to {\Bbb T}$ defined by the formula  
         \begin{equation}
z =\frac{\lambda-i}{\lambda+i}=: \omega (\lambda).
\label{eq:CV}\end{equation}
We transplant the measure $dm (\lambda)$ on  ${\Bbb R}$ to the measure $d\mu (z)$ on  ${\Bbb T}$ by the formula
\begin{equation}
 \mu (Y)=m (\omega^{-1} (Y))
\label{eq:CV1}\end{equation}
for all Borelian sets $Y\subset{\Bbb T}$. In particular, it follows from \e{eq:CV1} that
\begin{equation}
\int_{\Bbb T} z^{-n} d\mu(z)= \int_{\Bbb R} \Big( \frac{\lambda-i}{\lambda+i}\Big)^{-n} dm(\lambda),\q \forall n\in{\Bbb Z}.
\label{eq:CV2}\end{equation}

  \medskip
  
  {\bf 4.2.}
Let us now  express the right-hand side of   \e{eq:CV2} in terms of the Laguerre polynomials (see the book \cite{BE}, Chapter~10.12)  defined by the formula
\begin{equation}
{\sf L}_{n}^\alpha  (x)= n!^{-1} e^x  x^{-\alpha} d^n (e^{-x} x^{n +\alpha }) / dx^n, \q n=0,1, \ldots, \q \alpha\geq 0 .
\label{eq:K3xb}\end{equation}
Of course the polynomial ${\sf L}_{n} ^\alpha (x)$ has  degree $n$. We recall that  the polynomials ${\sf L}_{n} ^\alpha (x)$ are obtained by the orthogonalization of the monomials $1,x, x^2,\ldots$ with respect to the scalar product
\begin{equation}
{\pmb \la}f_{1}, f_{2}  {\pmb \ra}_{\alpha}=\int_{{\Bbb R}_{+}}  f_{1}(x) f_{2}(x) x^\alpha e^{-x} dx.
\label{eq:weight}\end{equation}
We need these polynomials for $\alpha=0$  and $\alpha=1$ only.
It follows from \e{eq:weight} that
 \begin{equation}
  \int_{{\Bbb R}_{+}}     {\sf L}_{n}^1 (x ) x^q e^{-   x}   d x= 0 \q{\rm if}\q  1\leq q\leq n.
\label{eq:La}\end{equation}

Recall the identity (see formula (10.12.32) in \cite{BE})
 \[
  \int_{{\Bbb R}_{+}}    {\sf L}_{n}^1 (x ) x e^{- \zeta x}  d x=  (n+1)   \frac{ (\zeta-1)^n}{ \zeta^{n+2}}, \q \Re \zeta> 0 ,
\]
whence
 \begin{equation}
\frac{d}{d\z}  \int_{{\Bbb R}_{+}}      {\sf L}_{n}^1 (x )  e^{- \zeta x} d x= - (n+1)   \frac{ (\zeta-1)^n}{ \zeta^{n+2}}  .
\label{eq:KL4x5}\end{equation}
  Since
\[
\frac{d}{d\z} \Big(\frac{ \zeta-1 }{ \zeta }\Big)^{n+1}  =   (n+1) \frac{ (\zeta-1)^n}{ \zeta^{n+2}},
\]
it follows from \e{eq:KL4x5} that, for some constant $c$,
 \begin{equation}
   \int_{{\Bbb R}_{+}}      {\sf L}_{n}^1 (x ) e^{- \zeta x}  d x= -  \Big(\frac{ \zeta-1 }{ \zeta }\Big)^{n+1} +c.
\label{eq:KL4x6}\end{equation}
Considering now the limit $\z\to + \infty$, we see that $c=1$. In particular, setting $\z= (1+i\lambda)/2$, we obtain the following result.

\begin{lemma}\label{Laguerre}
     For all $n=1,2,\ldots$, we have the representation
      \begin{equation}
      \Big( \frac{ \lambda - i }{ \lambda + i  }\Big)^{-n} = 1-2
        \int_{{\Bbb R}_{+}}      {\sf L}_{n-1}^1 (2x ) e^{-   x-i\lambda x}  d x .
\label{eq:Lag}\end{equation}
    \end{lemma}
    
 Putting in \e{eq:KL4x6} $\z=1$, we also see that
    \begin{equation}
  \int_{{\Bbb R}_{+}}     {\sf L}_{n}^1 (x ) e^{-   x}   d x= 1, \q n=0, 1,2,\ldots.
\label{eq:La1}\end{equation}

  Let $d m (\lambda)$ be an arbitrary complex measure such that $|m|({\Bbb R}) <\infty$. Integrating equality \e{eq:Lag} and
  using the Fubini theorem to interchange the order of integrations in $\lambda$ and $x$, we see that
   \begin{equation}
  \int_{\Bbb R}    \Big( \frac{ \lambda - i }{ \lambda + i  }\Big)^{-n} dm(\lambda) = m ({\Bbb R}) -2
        \int_{{\Bbb R}_{+}}     {\sf L}_{n-1}^1 (2x ) e^{-   x } \big(  \int_{\Bbb R} e^{-i\lambda x}   dm(\lambda)\big) d x .
\label{eq:Lag1}\end{equation}

  \medskip
  
  {\bf 4.3.}
    To give an idea of the proof of Theorem~\ref{RBT}, we first consider  the particular case $p=0$ when $dM(\lambda)= d m(\lambda)$ according to \e{eq:WHs2},  and  the measure $| M |(\Bbb R)=   |m|(\Bbb R) <\infty$. Let us also suppose that 
    the function \e{eq:WHs} equals zero for all $x>0$. 
       We proceed from equality \e{eq:Lag1}.
  By our assumption, the integral over $\lambda$ in the right-hand side   is zero for $x>0$. Recall that the measure $d\mu(z)$ on $\Bbb T$ is defined by equations   \e{eq:CV} and  \e{eq:CV1}; in particular, $\mu(\{1\})=0$.  Putting now   relations \e{eq:CV2}     and \e{eq:Lag1} together,  we see that
  \[
\int_{\Bbb T} z^{-n} d\mu(z)= m ({\Bbb R}), \q n=1,2, \ldots.
\]
Therefore Corollary~\ref{brothers1} implies that the measures $ d\mu (z) $ and hence  $dM(\lambda)$ are  absolutely continuous.

  \medskip
  
  {\bf 4.4.}
  Let us pass to the general case.     Observe first that replacing $dM(\lambda)$ by the measure  $e^{i a\lambda}dM(\lambda)$, we can suppose that $a=0$ in the conditions of Theorem~\ref{RBT}.
  Let the function  $\sigma(x)$ be defined by formula \e{eq:WHs}; then $\sigma\in L^2 ({\Bbb R}_{+})$. We set
  \[
  \wh{\sigma} (\lambda)= \int_{{\Bbb R}_{+}} \sigma(x) e^{ix\lambda} dx
  \]
  and introduce the measure
      \[
  d\wt{M}(\lambda)=  d{M}(\lambda) -\wh{\sigma} (\lambda) d \lambda.
  \]
  By the Parseval theorem,   $\wh{\sigma}\in L^2 ({\Bbb R})$ so that the measure
$  d\wt{M}(\lambda) $ satisfies condition \e{eq:Schs}  (for $p>1/2$) and
  \begin{equation}
\wt\sigma(x): = \frac{1}{2\pi} \int_{\Bbb R}  e^{-ix \lambda} d\wt{M}(\lambda ) =  \frac{1}{2\pi} \int_{\Bbb R} e^{-ix \lambda} d M (\lambda ) -\sigma(x)  =0
 \label{eq:WHsz}\end{equation}
 for $x>0$.   We now put
      \begin{equation}
        dm(\lambda )= (1+\lambda^2)^{-p} d\wt{M}(\lambda ) 
     \label{eq:Mm}\end{equation}
         and
           \begin{equation}
 \eta(x)  =
\frac{1}{2\pi} \int_{\Bbb R} e^{-ix \lambda}   dm(\lambda ).
 \label{eq:WHs1}\end{equation}
 Then $|m|({\Bbb R})<\infty$, $\eta\in L^\infty ({\Bbb R})\subset {\cal S}'$ and   (without loss of generality we suppose that 
 $p$ is an integer)
      \[
(-\frac{d^2}{d x^2} +1)^p \eta(x)  = \frac{1}{2\pi} 
 \int_{\Bbb R} e^{-ix \lambda} (1+\lambda^2)^p  dm(\lambda )=   \wt\sigma(x).
 \]
 Therefore  equation \e{eq:WHsz} for $x>0$ yields the differential equation
  \begin{equation}
(-\frac{d^2}{d x^2} +1)^p \eta(x)   =0, \q x>0,
 \label{eq:S}\end{equation}
 for the function \e{eq:WHs1}.
 Its solutions  are linear combinations of classical functions $x^q e^{-x}$ and $x^q e^{x}$ where $q=0,1,\ldots, p-1$. Since the functions $x^q e^{x}$  do not belong to the class $ {\cal S}'$, the general solution  of equation
 \e{eq:S} in $ {\cal S}'$ is given by the formula
    \begin{equation}
 \eta (x)= \sum_{q=0}^{p-1} c_{q} x^q e^{-x} ,\q x>0,
\label{eq:S1}\end{equation} 
where $c_{q}$, $q=0,1,\ldots, p-1$,  are arbitrary complex numbers.
  
    Let us again use relation \e{eq:Lag1} which in view of \e{eq:WHs1} and \e{eq:S1} can be written as
   \[
  \int_{\Bbb R}    \Big( \frac{ \lambda - i }{ \lambda + i  }\Big)^{-n} dm(\lambda) = m ({\Bbb R}) -4\pi
 \sum_{q=0}^{p-1} c_{q}       \int_{{\Bbb R}_{+}}     {\sf L}_{n-1}^1 (2x ) x^q e^{-   2x }   d x .
\]
 Suppose that $n\geq p$. Then it follows from \e{eq:La} that all integrals in the right-hand side, except that for $q=0$, equal to zero. Therefore according to \e{eq:La1}  we have
   \[
  \int_{\Bbb R}    \Big( \frac{ \lambda - i }{ \lambda + i  }\Big)^{-n} dm(\lambda) = m ({\Bbb R}) -2\pi
  c_{0} , \q n\geq p      .
\]
Let again the measure $d\mu (z)$ on  ${\Bbb T}$ be defined by  formulas \e{eq:CV}, \e{eq:CV1}; in particular, $\mu(\{1\}) = 0$.
Using \e{eq:CV2}  we now find that  condition \e{eq:br} is satisfied for $N=p -1$ and $b =  m ({\Bbb R}) -2\pi
  c_{0}$. Therefore, by Corollary~\ref{brothers1}, the  measures $dm(\lambda)$ and hence $dM(\lambda)$ are absolutely continuous. This concludes the proof of Theorem~\ref{RBT}. $\q\Box$

 \begin{remark}\label{Mm}
Of course the choice \e{eq:CV} of the mapping $\omega \colon{\Bbb T}\to {\Bbb R}$ is not unique. If however
$\omega (\lambda) =\frac{\lambda-i a}{\lambda+ia}$
  where $a>0$, then the proof above requires that the definition \e{eq:Mm}  of the measure $dm(\lambda)$ be also changed to $    dm(\lambda )= (a^2+\lambda^2)^{-p} d\wt{M}(\lambda ) $.
    \end{remark}
    
  \section{Closable Wiener-Hopf quadratic forms}  
%

In the  first subsection we   prove Theorem~\ref{T1}.  The next two subsections consist of its discussion and of an example of a highly singular but closable Wiener-Hopf quadratic form.

\medskip

{\bf 5.1.}
 The forms $w[f,f]$ and ${\bf w}[f,f]$ are defined by the same expression \e{eq:QFx}, but
 \[
 {\cal D}[w]=C_{0}^\infty ({\Bbb R}_{+} )\subset C_{0}^\infty ({\Bbb R}  )=  {\cal D}[{\bf w}].
 \]
 Therefore if ${\bf w}[f,f]$ is closable, then the same is true for $w[f,f]$.
 Thus  the ``$\,$if$\,$" part of Theorem~\ref{T1} is a consequence of the ``$\,$if$\,$"  part of Theorem~\ref{T1B}.
 
 The  proof of the converse statement follows   the same scheme as in Theorem~\ref{T1B}, but it is   more involved analytically since   Theorem~\ref{RBT} (the continuous analogue of the Riesz Brothers theorem) is required now.
 Below we avoid repeating the   arguments already used by the proof of Theorem~\ref{T1B}.
 
 We   suppose that estimate \e{eq:T1} is true with $\gamma=1$ 
 and denote by $d {\sf M}(\lambda )$ the measure satisfying equation \e{eq:WH}.
      Let us define   an   operator $A \colon L^2 ({\Bbb R}_{+})\to    L^2 (\Bbb R; d {\sf M} )$  on domain $\cal D (A)=  C_{0}^\infty ({\Bbb R}_{+})$ by the formula
   \[
(Af) (\lambda)=\frac{1}{\sqrt{2\pi}}\int_{{\Bbb R}_{+}} e^{ix\lambda} f(x) d x ,\q f\in C_{0}^\infty ({\Bbb R}_{+}) .
\]
 In view of     \e{eq:WH} the form       \e{eq:QFq}   can be written (cf. the identity  \e{eq:A2}) as
   \[
w[f,f] =  \| Af\|^2_{ L^2 (\Bbb R; d{\sf M})} , \q f\in C_{0}^\infty ({\Bbb R}_{+}).
\]
This  yields the following result (cf. Lemma~\ref{de}).
 
    \begin{lemma}\label{deW}
 The form $ w [f,f]$ defined   on $C_{0}^\infty ({\Bbb R}_{+})$   is closable in the space $L^2 ({\Bbb R}_{+} )$  if and only if the operator $A \colon L^2 ({\Bbb R}_{+} )\to    L^2 (\Bbb R; d{\sf M})$ defined   on the  domain $\cal D (A)=C_{0}^\infty ({\Bbb R}_{+})$  is closable.  
    \end{lemma}
    
  The function (distribution)  $u_{*} $ is again defined by formula \e{eq:A1}, but instead of Definition~\ref{def} we now accept
  
     \begin{definition}\label{defW}
 The set   ${\cal D}_{*}\subset   L^2 (\Bbb R; d {\sf M})$ consists  of all  $u \in  L^2 (\Bbb R; d {\sf M})$ such that   $u_{*} \in  L^2 ({\Bbb R}_{+} )$.
    \end{definition}
    
    With this definition of   ${\cal D}_{*} $, 
 Lemma~\ref{LTM} remains unchanged, and the following result plays the role  of  Lemma~\ref{adj}.

     \begin{lemma}\label{adjW}
 The operator $A $ and the form $w[f,f]$ are  closable if and only if the relation \e{eq:D} holds.
    \end{lemma}

     The central point is to check that the measures $ u(\lambda) d {\sf M} (\lambda)$ are absolutely continuous for all $u  \in     {\cal D}_{*} $. For such $u$, we have $u_{*} \in  L^2 ({\Bbb R}_{+})$, and so this fact follows from  Theorem~\ref{RBT} applied to  the measure $ d M (\lambda)= u(\lambda) d {\sf M} (\lambda)$.

  Now   we only have to repeat the arguments used for the proof
      of  Theorem~\ref{T1B}. 
If the form $ w [f,f]$  is closable, then by Lemma~\ref{adjW} the condition  \e{eq:D} is satisfied. 
Since   the measures $    u(\lambda) d {\sf M} (\lambda)$are absolutely continuous,  the same is true for the measure $  d {\sf M} (\lambda)$  because   Lemma~\ref{ac} remains  obviously true.       $\Box$

    Comparing Theorems~\ref{T1}  and \ref{T1B}  we see that the form $w[f,f]$  is closable if and only if
    the form ${\bf w}[f,f]$  is closable.
    
    \medskip

{\bf 5.2.}
Let us state a consequence of Theorem~\ref{T1} in terms of the integral kernel $w(x)$   of the form \e{eq:QFq}.

  \begin{proposition}\label{T1m} 
Suppose that the condition 
 \e{eq:T1} is satisfied. 
If  the form $w[f,f] $  is closable  in the space $L^2 ({\Bbb R}_{+})$, then
 \begin{equation}
  w(x)= (1+D^{2})^{q} v(x)
  \label{eq:Ker}\end{equation}
   for some number $q \in  {\Bbb Z}_{+}$ and a function $v\in L^\infty ({\Bbb R} )$ such that $v(x)\to 0$ as $|x |\to \infty$.
   If condition \e{eq:Ker} is satisfied for some  $q \in  {\Bbb Z}_{+}$ and   $v\in L^2 ({\Bbb R} )$, then 
   the form $w[f,f] $  is closable  in the space $L^2 ({\Bbb R}_{+})$.
    \end{proposition}
    
     \begin{pf}
   If  the form $w[f,f] $  is closable, then, by   Theorem~\ref{T1}, we have
    \begin{equation}
  w(x)= \frac{1}{2\pi} \int_{\Bbb R} e^{-ix\lambda} (1+\lambda^{2})^{p} \psi (\lambda)d\lambda
  \label{eq:Ker1}\end{equation}
  for some number $p$ and a function   $\psi\in L^1 ({\Bbb R} )$. Therefore representation \e{eq:Ker} holds with any integer  $q\geq p$ and the function $v =(2\pi)^{-1/2} \Phi \psi_{1}$ where $\psi_{1}(\lambda)=(1+\lambda^{2})^{p-q} \psi (\lambda)$. 
  The function $v (x)$ is 
  bounded and tends to zero as $|x |\to \infty$ by the Riemann-Lebesgue lemma.
    
     On the other hand,  it follows from \e{eq:Ker} where $v\in L^2 ({\Bbb R} )$ that  representation \e{eq:WH} is true with 
the measure
        \[
  d{\sf M}(\lambda)= (2\pi)^{1/2}  (1+\lambda^{2})^{q} (\Phi^{* }v) (\lambda)  d \lambda.
 \]
 Since  $ \Phi^{* }v \in L^2 ({\Bbb R}_{+})$, condition  \e{eq:Sch} on this measure is satisfied for $p>q+1/4$.
        \end{pf}
        
        There is an obvious gap between necessary  and sufficient conditions in Proposition~\ref{T1m}. The conditions 
        $w\in L^\infty ({\Bbb R} )$ and $w(x) \to 0$ as $|x |\to \infty$ do not imply that the Fourier transform $\wh{w}
        \in L^1_{\rm loc} ({\Bbb R} )$. So these conditions are not sufficient.  On the other hand, the condition \e{eq:Ker} with $v\in L^2 ({\Bbb R} )$ implies that  $\wh{w}
        \in L^2_{\rm loc} ({\Bbb R} )$ which is stronger than $\wh{w}
        \in L^1_{\rm loc} ({\Bbb R} )$.
        
        For Toeplitz operators, the discrete analogue of Proposition~\ref{T1m} is discussed in \cite{YaT} in more details.
        
           \medskip

{\bf 5.3.}
We emphasize that a relatively difficult part of Theorem~~\ref{T1} is its ``$\,$only if$\,$" assertion. However the ``$\,$if$\,$" part also gives interesting examples of highly singular but closable Wiener-Hopf quadratic forms.

 \begin{example}\label{exX} 
 Set
 \[
 w(x)=    e^{-\pi ia /2} x_{+}^{-a} +  e^{\pi ia /2} x_{-}^{-a} ,\q \forall a> 0,
 \]
 where the distribution $x_{\pm}^{-a}$ is  standardly defined (see, for example, \cite{GUEVI+}) by the analytic continuation in the parameter $a$ of the   function $| x |^{-a}\1_{{\Bbb R}_\pm} (x)$ which belongs to $L^1_{\rm loc}$ for $a\in ( 0,1)$.
 Note   that $w(x)=   w_{a}(x)$ is well defined as a distribution in the Schwartz space ${\cal S}'$ for all values of $a>0$ including integer values. For example, for $a=2$ we have
 \[
 \la w_2 ,\theta\ra= - \int_{{\Bbb R}_{+}}x^{-2}\big(\, \ov{\theta (x)} +\ov{\theta (-x)}\,-2 \, \ov{\theta (0)} \, \big) dx -\pi i \, \ov{\theta' (0)}   .
 \]
 Calculating the Fourier transform of the function $w(x)$, we see that relation \e{eq:WH}  is satisfied with the absolutely continuous measure
\begin{equation}
d{\sf M} (\lambda)= 2\pi  \Gamma (a )^{-1}\1_{{\Bbb R}_{+}} (\lambda) \lambda^{a-1} d\lambda 
\label{eq:MM}\end{equation}
      where $\Gamma (\cdot)$ is the gamma function.
Thus according to Theorem~~\ref{R-H} the Wiener-Hopf quadratic form $w[f,f]$ is non-negative. It follows from  
Theorem~~\ref{B-H} that the corresponding operator $W$ is unbounded unless $a=1$. Finally,
  Theorem~~\ref{T1} implies that the form $w[f,f]$ is closable for all values of $a>0$.
    \end{example}

\section{Discussion}  

   In Subsection~6.1, we compare our results with the similar statements for Toeplitz  operators obtained in \cite{YaT}.  
    In Subsection~6.2, we discuss
      the corresponding assertions for Hankel   operators (realized as integral operators in the space $L^{2} ({\Bbb R}_{+})$) obtained in \cite{Yf1a}. 
  
  \medskip

{\bf 6.1.} 
  As was already noted, Wiener-Hopf  $W$ and Toeplitz $T$ (see definition \e{eq:HDT}) operators    can be considered as  continuous and discrete versions of the same object. All the results stated in Section~1, have their counterparts in the theory of Toeplitz operators. Similarly to  \e{eq:QFq}, the Toeplitz quadratic form is given by the equation
   \begin{equation}
t[g,g] =\sum_{n,m\geq 0} t_{n-m} g_{m}\bar{g}_{n}.   
 \label{eq:QFqt}\end{equation} 
 Here sequences $g=\{g_{n}\}_{n\in {\Bbb Z}_{+}} $ have only a finite number of non-zero components. The set of such $g$ is denoted $\cal D $. It is dense in $\ell^2({\Bbb Z}_{+})  $. A priori there are no restrictions on the sequence $t=\{t_{n}\}_{n\in {\Bbb Z} } $ in  \e{eq:QFqt} except that   $t_{n}= \ov{t_{-n}}$ so that the quadratic form $t[g,g] $ is real.

  In the theory of Toeplitz   operators, the F.~Riesz-Herglotz theorem   plays exactly the same role as the
  Bochner-Schwartz theorem      plays in the theory of Wiener-Hopf operators. It states that   $t[g,g]\geq 0$ for all 
  $g\in\cal D $ if and only if
  \begin{equation}
  t_{n} =\int_{\Bbb T} z^{-n} d\mu (z)  , \q \forall n\in{\Bbb Z},
 \label{eq:WHt}\end{equation}
 for some 
  non-negative  measure $d\mu(z)$ on the unit circle $\Bbb T$.   
  
  The role of Theorem~\ref{T1} is played by the following assertion.
  
   \begin{theorem}\label{T1t}\cite[Theorem 1.3]{YaT}   
    Let the form $t[g,g]$ be given by formula \e{eq:QFq} on elements $g\in {\cal D}$, and let the condition 
    \[
t[g,g] \geq \gamma \| g\|^2  ,\q g\in{\cal D}, \q \| g\| =\| g\|_{\ell^2({\Bbb Z}_{+})},
 \]
  be satisfied for some   $\gamma\in \Bbb R$. 
   Then the form $t[g,g] $ is closable in the space $\ell^2 ({\Bbb Z}_{+})$ if and only if the measure $d\mu  (z)$ in   the equations \e{eq:WHt} is absolutely continuous.
    \end{theorem}
      
   Morally,    Theorems~\ref{T1} and \ref{T1t} are of course quite similar, but  mathematically they are not completely equivalent.
   Let us explain the link and the difference between these results. Set
   \begin{equation}
   (U g)(x)=\sqrt{2}\sum_{n \in{\Bbb Z}_{+}}g_{n}{\sf L}^{0}_{n} (2x) e^{-x}, \q g =\{g_{n}\}_{n\in {\Bbb Z}_{+}}\in {\cal D},
  \label{eq:WH1}\end{equation}
  where the Laguerre polynomials ${\sf L}^{0}_{n} (x)$ are defined by formula \e{eq:K3xb}.
   Since the functions $\sqrt{2} {\sf L}^{0}_{n} (2x) e^{-x }$, $n\in{\Bbb Z}_{+}$, form an orthonormal basis in $ L^{2} ({\Bbb R}_{+}) $,   operator \e{eq:WH1}  extends to
     the unitary mapping $U: \ell^{2} ({\Bbb Z}_{+})\to L^{2} ({\Bbb R}_{+})$.
      Suppose that the function (distribution) $w(x)$ is given by formula  \e{eq:WH}. 
          Substituting  \e{eq:WH} and \e{eq:WH1} into expression  \e{eq:QFq} for the Wiener-Hopf quadratic form, we find that
       \begin{equation}
   w[U g, Ug]=\frac{1}{\pi}\sum_{n,m \in{\Bbb Z}_{+}}g_{m}\ov{g_n}\int_{  {\Bbb R}_{+}}\int_{  {\Bbb R}_{+}} dxdy{\sf L}^{0}_{n} ( 2 x) {\sf L}^{0}_{m} (2 y) e^{-x-y } \Big(\int_{\Bbb R} e^{-i \lambda(x-y)} d{\sf M} (\lambda)\Big) .
  \label{eq:WH2}\end{equation}
  Interchanging here the order of integrations in $x,y$ and $\lambda$ and using the formula (see  (10.12.32) in \cite{BE})
 \[
 2 \int_{{\Bbb R}_{+}}     {\sf L}_{n}^0 ( 2 x ) e^{- 2\zeta x}  d x=      \frac{ (\zeta-1)^n}{ \zeta^{n+1}}, \q \Re \zeta> 0 ,
\]
where $\z= (1+i\lambda)/2$
for the integrals over $x$ and $y$, we see that
     \begin{equation}
   w[U g, Ug]= \frac{1}{\pi}\sum_{n,m \in{\Bbb Z}_{+}}g_{m}\ov{g_n} \int_{\Bbb R} (\lambda^{2}+1)^{-1}\Big(\frac{\lambda -i}
   {\lambda +i}\Big)^{n-m}d{\sf M} (\lambda)  .
  \label{eq:WH3}\end{equation}
  Now it follows from \e{eq:QFqt} that
   \begin{equation}
  w[U g, Ug]=t[g,g]
    \label{eq:WH4}\end{equation}
    provided
      \begin{equation}
   t_{n}= \frac{1}{\pi} \int_{\Bbb R} (\lambda^{2}+1)^{-1}\Big(\frac{\lambda -i}
   {\lambda +i}\Big)^{-n}d{\sf M} (\lambda)  .
  \label{eq:WH3a}\end{equation}
  Making here the change of variables \e{eq:CV}, we see that this expression coincides with \e{eq:WHt} if
 the measure $d\mu (z)$ on the unit circle is given by the formula
    \begin{equation}
 \mu (Y)= \frac{1}{\pi} \int_{\omega^{-1} (Y)} (\lambda^{2}+1)^{-1} d{\sf M} (\lambda) .
\label{eq:CV1T}\end{equation}
The arguments above do not require that the measures be  non-negative; they may be even complex.
Note however that the measure $d\mu (z)$ is non-negative (semibounded) if and only if the measure ${\sf M} (\lambda)$ is non-negative (semibounded).
  
     According to formula \e{eq:WH4} one might think that the form $w[f,f]$ defined on $C_{0}^{\infty} ({\Bbb R}_{+})$ is closable if and only if this is true for the form $t[g,g]$ defined on ${\cal D}$ provided $w(x)$ and $t_{n}$ are given by formulas \e{eq:WH} and \e{eq:WH3a}, respectively.  However this is not completely true. First of all, we note that $U{\cal D}\neq C_{0}^{\infty} ({\Bbb R}_{+})$ so that the domains of $t[g,g]$ and $w[f,f]$ are not linked by the mapping $U$. Second, it follows from \e{eq:CV1T} that    
     $ |\mu |({\Bbb T})<\infty$ and formula \e{eq:WH3a} makes sense if condition  \e{eq:Schs} is satisfied for $p=1$, but not for larger $p$. Finally, the passage from \e{eq:WH2} to \e{eq:WH3} can be justified by the Fubini theorem, but it requires that  $ |{\sf M} |({\Bbb R})<\infty$,  which is true only if $p=0$ in \e{eq:Schs}.  Since $p$ is arbitrary  in Theorem~\ref{T1}, it is more general than Theorem~\ref{T1t}.

 \medskip

{\bf 6.2.}
There is a certain parallelism between theories of Toeplitz and Hankel operators. This is true both for the discrete (in the space $\ell^{2} ({\Bbb Z}_{+})$) and for the continuous (in the space $L^{2} ({\Bbb R}_{+})$) realizations of these operators. Here we discuss the continuous realizations; see \cite{YaT}, for a discussion of the discrete realizations. For example, the criteria of  boundedness
of Wiener-Hopf and of Hankel operators due to Toeplitz  (see Theorem~\ref{B-H}) and to Nehari \cite{Nehari}, respectively, look formally similar. In the semibounded case, the study of  
Wiener-Hopf  quadratic forms relies on  the Bochner theorem   on the Fourier transform of  functions of positive type while the study of   Hankel quadratic forms relies on the  Bernstein theorem on  exponentially convex functions. These results play the role of the trigonometric and power moment problems, respectively, for the discrete realizations.

To be more precise, we use the generalization by L.~Schwartz of the Bochner theorem to distributions (see Theorem~\ref{R-H}). Similarly,
for applications to Hankel operators, we need a generalization to distributions of the Bernstein theorem. Let us state it here. For $f,g\in C_{0}^\infty({\Bbb R}_{+})$, we set
\[
 (f\star g)(t)=\int_{0}^t f(s) g(t-s) ds.
 \]
 Obviously, $f \star g\in C_{0}^\infty({\Bbb R}_{+})$.

\begin{theorem}\label{Bern}\cite[Theorem~5.1]{Yafaev3}
Let $h\in C_{0}^\infty({\Bbb R}_{+})'$ and
\begin{equation}
h[f,f]:= {\pmb \la}h, \bar{f}\star f {\pmb \ra}\geq 0
\label{eq:BernH}\end{equation}
for all $f\in C_{0}^\infty({\Bbb R}_{+})$. Then there exists a non-negative measure $d{\sf M} (\lambda) $ on $\Bbb R$ such that
\begin{equation}
h(t)= \frac{1}{2\pi}\int_{\Bbb R} e^{-t\lambda} d {\sf M}  (\lambda)  
\label{eq:Bern1}\end{equation}
where  the integral converges for all $t> 0$.
 \end{theorem}
 
 We emphasize that the measure $d{\sf M} (\lambda) $ may grow almost exponentially as $\lambda\to +\infty$ and it tends to zero super-exponentially as $\lambda\to -\infty$, that is,
 \[
 \int_{{\Bbb R}_{+}} e^{-t\lambda} d{\sf M} (\lambda)<\infty \q {\rm and}\q   \int_{{\Bbb R}_{+}} e^{t\lambda} d{\sf M} (-\lambda) <\infty
\]
  for an arbitrary small $t>0$ and  for an arbitrary large $t>0$, respectively. In the theory of Hankel operators, representation \e{eq:Bern1} plays the same role as representation \e{eq:WH} plays in the theory of Wiener-Hopf operators.
 
   Theorem~\ref{Bern} shows that  the positivity of ${\pmb\la}h,\bar{f} \star f {\pmb\ra}$ imposes very strong conditions on $h(t)$. They are stated in the   assertion below.

 \begin{corollary}\label{BernN}
 Under the assumptions of Theorem~\ref{Bern},
the function  $h\in C^\infty ({\Bbb R}_{+})$. Moreover, it  admits the analytic continuation in the right-half plane $\Re t> 0$ and   is uniformly bounded in every strip $\Re t \in (t_{1},t_{2})$ where $0< t_{1} <t_{2} <\infty$.
 \end{corollary}
 
 We recall  also the following result.
 
    \begin{theorem}\label{Carl}\cite[Theorem~3.10]{Yf1a}
     Let   $h(t)$ be given by formula \e{eq:Bern1} where ${\sf M} ((-\infty,0])=0$. Then the form \e{eq:BernH} defined on $f\in C_{0}^{\infty}({\Bbb R}_{+})$ is closable.
        \end{theorem}
        
       If $h(t)$ is given by   \e{eq:Bern1}, then 
       \[
       h(t)\geq  {\sf M} ((-\infty,0]), \q \forall t>0.
       \]
       Therefore putting together Theorems~\ref{Bern} and \ref{Carl}, we obtain  a simple sufficient condition for a Hankel quadratic form to be closable.

   \begin{theorem}\label{Hamb1} 
     Let   assumption \e{eq:BernH} be satisfied. Then the form $h[f,f]$ defined on $f\in C_{0}^{\infty}({\Bbb R}_{+})$ is closable if
     \begin{equation}
\liminf_{t \to\infty}  h(t)=0.
\label{eq:Carl}\end{equation}
In thi                                                s case, there exists a unique non-negative operator $H$ corresponding $($see Subsection~$2.1)$ to the form $h[f,f]$.
        \end{theorem}
        
          \begin{remark}\label{inf} 
  If condition  \e{eq:Carl} is satisfied, then ${\sf M} ((-\infty,0])=0$, and therefore representation \e{eq:Bern1} implies that $h(t)\to 0$   monotonically as $t\to \infty$. Actually, such functions $h(t)$ are called (see  \cite{AKH}, Chapter~5, \S 5)                                completely monotonic.          \end{remark}
        
          \begin{example}\label{Carl1} 
     Let   $h(t)=t^{-a}$ where $a>0$. Then $h[f,f]\geq 0$, and the form $h[f,f]$   is closable.  
             \end{example}
             
             Note that in this example the representation \e{eq:Bern1} is satisfied with the measure $  d{\sf M}(\lambda)$
           given by formula \e{eq:MM}. The corresponding
        Hankel operators $H$  are unbounded  unless $a = 1$; see \cite{Yf1a} for details. We emphasize  that the singularity of $h(t)$ as $t\to 0$ may be arbitrary strong. On the other hand,  according to Example~\ref{ex}  the form $h[f,f]$ is not closable if $ h(t)=1$. So the sufficient condition \e{eq:Carl} is rather close to necessary.

       As observed at the end of Section~5, for Wiener-Hopf quadratic forms, the condition $w(x)\to 0$ as $|x|\to \infty$ does not imply that $w[f,f]$ is closable. According to Theorem~\ref{Hamb1} the situation is different for  Hankel quadratic forms.


\end{document}